\documentclass[11pt]{article}
\PassOptionsToPackage{obeyspaces}{url}
\usepackage{amsfonts, amsmath, amssymb}

\usepackage[hidelinks]{hyperref}
\usepackage{bookmark}
\bookmarksetup{open,numbered,depth=5} 

\usepackage{amsthm}

\usepackage{breakurl}
\usepackage{tikz}

\usepackage{etoolbox} 
\patchcmd{\thebibliography}{\leftmargin\labelwidth}{\leftmargin\labelwidth\addtolength\itemsep{-0.1\baselineskip}}{}{}

\usepackage{mathtools}

\author{Boris Bukh\thanks{Department of Mathematical Sciences, Carnegie Mellon University, Pittsburgh, PA 15213, USA\@. Supported in part by U.S.\ taxpayers through NSF grant DMS-2154063.}\and
R. Amzi Jeffs\thanks{Department of Mathematical Sciences, Carnegie Mellon University, Pittsburgh, PA 15213, USA\@. Supported by the National Science Foundation through Award No. 2103206.}}

\title{Distances between realizations of order types}
\date{September 2023}

\usepackage[nameinlink]{cleveref}

\newtheorem{theorem}{Theorem}
\newtheorem{lemma}[theorem]{Lemma}
\newtheorem{prop}[theorem]{Proposition}

\newcommand*{\R}{\mathbb{R}}                                     
\newcommand*{\eqdef}{\stackrel{\mbox{\normalfont\tiny def}}{=}}  
\newcommand*{\veps}{\varepsilon}                                 
\newcommand*{\Type}{\mathcal{T}}                                 
\DeclarePairedDelimiter\abs{\lvert}{\rvert}                      
\DeclareMathOperator{\orient}{orient}                            
\DeclareMathOperator{\OT}{OT}                                    
\DeclareMathOperator{\dist}{dist}                                
\DeclareMathOperator{\diam}{diam}                                
\DeclareMathOperator{\vol}{vol}                                
\DeclareMathOperator{\conv}{conv}                                
\DeclareMathOperator{\diag}{diag}                                
\DeclareMathOperator{\SO}{SO}                                    

\begin{document}

\maketitle

\begin{abstract}
  Any $n$-tuple of points in the plane can be moved to any other $n$-tuple by a continuous motion
  with at most $\binom{n}{3}$ intermediate changes of the order type.
  Even for tuples with the same order type, the cubic bound is sharp: there exist pairs of $n$-tuples
  of the same order type requiring $c\binom{n}{3}$ intermediate changes. 
\end{abstract}


\paragraph{Order types.}
A collection of $d+1$ points $p_0,\dotsc,p_d$ in $\R^d$ can be \emph{oriented} in one of the three ways according
to the sign of $\det \left[\begin{smallmatrix}1&\dotsb&1\\p_0&\dotsb&p_d\end{smallmatrix}\right]$.
We call this sign \emph{orientation} of $(p_1,\dotsc,p_{d+1})$ and denote it by $\orient(p_0,\dotsc,p_d)$. In the plane,
the three possible orientations correspond to the points lying in clockwise, counterclockwise and collinear positions.

\parshape=13 0cm\hsize 0cm.76\hsize 0cm.76\hsize 0cm.76\hsize 0cm.76\hsize 0cm.76\hsize 0cm.76\hsize 0cm\hsize 0cm\hsize 0cm\hsize 0cm\hsize 0cm\hsize 0cm \hsize
The totality of all $\binom{n}{d+1}$ orientations of points in a tuple $P=(p_1,\dotsc,p_n)$ of points in $\R^d$ determine the \emph{order type}, which we denote
by $\OT_P$.
\vadjust{\hfill\smash{\raise -64pt\llap{%
\begin{tikzpicture}[pt/.style={below,execute at begin node=$\scriptscriptstyle, execute at end node=$}]
\begin{scope}
\fill (-0.5,0) circle (1.6pt) node[pt] {1};
\fill (-0.1,1.1) circle (1.6pt) node[pt] {2};
\fill (0.5,0.3) circle (1.6pt) node[pt] {3};
 \fill (0,0.5) circle (1.6pt) node[pt] {4};
\end{scope}
\draw[densely dotted](1,-0.2) -- +(0,1.5);
\node at (1,-0.6) {\tiny Two tuples of};
\node at (1,-0.86) {\tiny the same order type};
\begin{scope}[xshift=2cm]
\fill (0.1,0) circle (1.6pt) node[pt] {1};
\fill (-0.6,0.35) circle (1.6pt) node[pt] {2};
\fill (0.5,1) circle (1.6pt) node[pt] {3};
\fill (0.2,0.6) circle (1.6pt) node[pt] {4};
\end{scope}
\end{tikzpicture}\quad}}}
For example, the two tuples of four points in the figure on the right have the same order type. It is tempting to
regard tuples with the same order type as similar. In particular, it was conjectured by Ringel~\cite{ringel}
that every two tuples of the same order type can be transformed into one another via a continuous motion,
without breaking the order type. Surprisingly, White \cite{white} disproved the conjecture by constructing an order type
with disconnected realization space. An even more surprising result was proved by Mn\"ev~\cite{mnev_phd}, which became known as the Mn\"ev's universality theorem.
It implies existence of order types whose realization spaces are homotopy equivalent to any prescribed simplicial complex, for instance.
The proof of Mn\"ev's theorem has since been simplified, see \cite{richter_gebert} for example.

\paragraph{Distance.} An $n$-tuple $P$ is a \emph{realization} of the order type $\Type\in \{-1,0,+1\}^{\binom{[n]}{d+1}}$ if $\OT_P=\Type$.
A tuple of points is in \emph{general position} if no $d+1$ its points lie on a common hyperplane. Similarly, we say that an
order type $\Type$ is in general position if $\Type\in \{-1,+1\}^{\binom{[n]}{d+1}}$.

In this paper, we seek ways to continuously transform one tuple in general position into another
with as few changes of the order type as possible. Formally,
a continuous function $m\colon [0,1]\to (\R^d)^n$ is a \emph{motion}
between tuples $P$ and $P'$ if $m(0)=P$, $m(1)=P'$, there are at most finite number of times $t$
for which $m(t)$ is not in general position, and at each of those times just there a unique set of $d+1$ points
lying on a common hyperplane. The number of such times is called the \emph{cost} of $m$.
Let $\dist(P,P')$ be the cost of the least expensive motion from $P$ to~$P'$. We call tuples $P$ and $P'$ \emph{isotopic}
if $\dist(P,P')=0$.

\paragraph{Our results.} It is clear that $\dist(P,P')$ is bounded from below by the Hamming distance
between $\OT_P$ and $\OT_{P'}$. In particular, if $P'$ is a mirror image of $P$, then $\OT_P=-\OT_{P'}$ and so
$\dist(P,P')\geq \nobreak \binom{n}{d+1}$ in this case. Our first result shows that, at least in the dimension $2$, this
is the largest distance between any two tuples.

\begin{theorem}\label{thm:genupper} Suppose $d\geq 2$. Then for any two $n$-tuples in general position $P,P'\in (\R^d)^n$ we have
\[
  \dist(P,P')\leq \frac{d}{2}\binom{n}{d+1}.
\]
\end{theorem}
We do not know if the fraction $\frac{d}{2}$ can be improved. In particular, we do not know if it can be replaced by an absolute
constant that is independent of~$d$.\smallskip

Though Ringel's conjecture is false, the following pair of results shows that its spirit is true --- it is easier to
transform one tuple into another if they have the same order type.
\begin{theorem}\label{thm:sametype}
If tuples $P,P'\in (\R^d)^n$ are in general position with the same order type and $d$ is even, then
  \[
    \dist(P,P')\leq \frac{d}{2}\binom{n}{d+1}-\frac{d}{2}.
  \]
\end{theorem}
We say that a tuple $P$ is \emph{$\alpha$-non-elongated} if the ratio between the largest and smallest
distances in $P$ is at most $\alpha\sqrt[d]{n}$.
\begin{theorem}\label{thm:nonelongated} Suppose $d\equiv 2\pmod 4$.
If tuples $P,P'\in (\R^d)^n$ are both $\alpha$-non-elongated and in general position with the same order type, then
  \[
    \dist(P,P')\leq \left(\frac{d}{2}-\beta\right)\binom{n}{d+1},
  \]
  for some positive constant $\beta=\beta(\alpha)>0$.  
\end{theorem}

On the other hand, there are tuples of the same order type that are hard to turn into one another.
\begin{theorem}\label{thm:blowup}
  There are arbitrarily large $n$-tuples $P,P'$ in $\R^2$
  of a same order type such that
  \[
  \dist(P,P')\geq 10^{-4}\binom{n}{3}.
  \]
\end{theorem}

\section*{Proof of~\texorpdfstring{\Cref{thm:blowup}}{Theorem 4}}
The requisite example is obtained from a counterexample to Ringel's conjecture
and the following proposition.
\begin{prop}\label{prop:blowup}
  Let $Q$ and $Q'$ be $r$-tuples of the same order type with $\dist(Q,Q')>0$. Then, for all sufficiently
  large $n$, there are $n$-tuples $P$ and $P'$ that that satisfy
  \[
    \dist(P,P')\geq \bigl(2r^{-3}+o(1)\bigr)\binom{n}{3}.
  \]
\end{prop}
The constant appearing in the statement of \Cref{thm:blowup} can be obtained both from the $13$-point example
of Tsukamoto \cite{tsukamoto} and from the earlier $14$-point example of Suvorov \cite{suvorov}. The rest of this section is devoted to the proof of~\Cref{prop:blowup}.\smallskip

We shall replace each point of $Q$ and of $Q'$ by a cloud of $m\eqdef n/r$ points.
The clouds will be sufficiently small so that whenever we select a point from each cloud,
the obtained tuples $\widetilde{Q}$ and $\widetilde{Q}'$ are isotopic to $Q$ and $Q'$ respectively.
It implies that at least one $(d+1)$-tuple changes orientation in any motion from $\widetilde{Q}$ to $\widetilde{Q}'$.
Double-counting such $(d+1)$-tuples according to the pair $(\widetilde{Q},\widetilde{Q}')$ that it belongs
would then yield
\[
  \dist(P,P')\geq 2m^{d+1}=2r^{-3}\binom{n}{3}+O(n^2).
\]
The factor of $2$ accounts for the fact that every $(d+1)$-tuple must change orientation an even number of times.\smallskip
  
The main difficulty is in selecting the clouds so that the resulting tuples $P$ and $P'$ still have the same order type.

\paragraph{Making clouds.}
A tuple $\widetilde{Q}=(\tilde{q}_1,\dotsc,\tilde{q}_r)$ is an $\veps$-perturbation of $Q=(q_1,\dotsc,q_r)$
if $\abs{\tilde{q}_i-q_i}\leq \nobreak \veps$. Pick $\veps>0$ sufficiently small so that every $\veps$\nobreakdash-perturbation of $Q$ is isotopic to $Q$,
and every $\veps$\nobreakdash-perturbation of $Q'$ is isotopic to~$Q'$.

Let $i\in [r]$ be arbitrary. Pick a unit vector $a_i$ that is not parallel to any vector of the form $q_j-q_i$ with $j\neq i$.
Let $L_i$ denote the directed line through $q_i$ in the direction $a_i$.
The $L_i$ partitions the rest of $Q$ into two sets $\{q_j\}_{j\in A}$ and $\{q_j\}_{j\in B}$, the left and right sides of the line.

Because $Q$ and $Q'$ have the same order type, the circular order of points around $q_i$ and $q_i'$ is the same, and so
we can find a unit vector $a_i'$ that such that the line $L_i'$ through $q_i'$ in the direction $a_i'$ induces the same partition
on $Q'$, namely $\{q_j'\}_{j\in A}$ and $\{q_j'\}_{j\in B}$.

Replace $q_i$ by the \emph{cloud} $Q_i$ consisting of points $q_{i,1},\dotsc,q_{i,m}$ where $q_{i,k}\eqdef q_i+k\delta a_i+k^2\delta^2 b_i$.
Similarly replace $q_i'$ by the points $q_{i,k}\eqdef q_i'+k\delta a_i'+k^2\delta^2 b_i'$ for $k\in [m]$. The number $\delta>0$ is
chosen sufficiently small to satisfy four conditions.
\begin{enumerate}
\item[(C1)] $\abs{q_{i,k}-q_i}\leq \veps$ holds for every $k$,
\item[(C1$'$)] $\abs{q_{i,k}'-q_i'}\leq \veps$ holds for every $k$,
\item[(C2)] if $q_{i,k}$ and $q_{i,k'}$ are two points of the same
cloud with $k<k'$ and $q_{j,\ell}$ is a point in another cloud, then $q_{j,\ell}$ is on the same side of the line $q_{i,k}q_{i,k'}$
as the side of $q_j$ relative to the line $L_i$,
\item[(C2$'$)] Same as (C2) but with the points of $Q_i$ and the line $L_i$ replaced by $Q_i'$ and the line $L_i'$.
\end{enumerate}
  
The requisite $n$-tuples are then $P\eqdef(q_{i,k})_{i\in [r],k\in [m]}$ and $P'\eqdef (q_{i,k}')_{i\in [r],k\in [m]}$.
Indeed, the lower bound on $\dist(P,P')$ follows from the fact that the cloud points are $\veps$\nobreakdash-perturbations of $Q$ and $Q'$
respectively, and the orientations of triples being the same is the consequence of the conditions (C2) and~(C2$'$).

\section*{Proof of \texorpdfstring{\Cref{thm:genupper}}{Theorem 1}}\label{sec:genupper}
We employ the usual arithmetic on tuples. For example, if $P=(p_1,\dotsc,p_n)$ is a tuple of points and $\lambda\in \R$,
then $\lambda P\eqdef (\lambda p_1,\dotsc,\lambda p_n)$. Similarly, if $M$ is a linear transformation, then $MP$ denotes the tuples obtained
by applying $M$ to each point in~$P$.

Given two tuples $P,P'$, a natural way to transform $P$ into $P'$ is via a \emph{linear motion}, whose parametrization is $tP'+(1-t)P'$ as time $t$
ranges from $0$ to~$1$. As we will see, doing so gives an upper bound of
$\dist(P,P')\leq d\binom{n}{d+1}$. The factor-of-two improvement in \Cref{thm:genupper} comes from first moving $P'$
into a suitable position.

The preparatory motion of $P'$ consists of two steps. First, using the fact that $P'$ is in general position, we perturb
points of $P'$ slightly. The precise conditions on the perturbation will be specified below (the determinants appearing in \Cref{lem:rj} must not vanish).
The second step, which consists of scaling the points of $P'$, we describe now.

Let $\diag(\lambda_1,\dotsc,\lambda_d)$ denote the diagonal $d$-by-$d$ matrix with entries $\lambda_1,\dotsc,\lambda_d$.
\begin{lemma}\label{lem:scale}
  Suppose that the numbers $\lambda_1,\dotsc,\lambda_d$ are nonzero, and the count of
  negative numbers among these is even. Let $D=\diag(\lambda_1,\dotsc,\lambda_d)$.
  Then $\dist(P,DP)=0$ for every tuple $P\in (\R^d)^n$.
\end{lemma}
\begin{proof}
  First suppose that all $\lambda$'s are positive. In this case,
  by continuously changing entries of the diagonal matrix we can turn the identity
  matrix into $D$. The resulting motion $t\mapsto D_t P$ evidently preserves the order type.

  Consider next the case $(\lambda_1,\dotsc,\lambda_d)=(-1,-1,+1,\dotsc,+1)$ (all but the first two entries are $+1$).
  In this case, let $D_t$ be the rotation around the origin in the $xy$-plane by the angle $\pi t$.
  Then $t\mapsto D_t P$ is the motion from $P$ to $P'$. We can similarly treat the case where two
  of the $\lambda$'s are equal to $-1$, and the rest are $+1$.

  The general case follows by combining these two cases. We first make all the signs positive, two at a time.
  Then we scale the magnitudes of the $\lambda$'s.
\end{proof}

Let $\lambda_1,\dotsc,\lambda_d$ be a sequence of non-zero real numbers satisfying the assumption of \Cref{lem:scale}, i.e.,
the number of negative elements is even. We postpone the actual choice.
Set $D\eqdef \diag(\lambda_1,\dotsc,\lambda_d)$, and let $P''=DP'$. In view of \Cref{lem:scale}, $\dist(P,P')=\dist(P,P'')$.

Consider what happens during the linear motion from $P$ to~$P''$. Since multiplication by a positive constant does not
change orientation, for $0\leq t<1$, the orientations of $(d+1$)-tuples in $(1-t)P+tP''$ is the same as in $P+\frac{t}{1-t}P''$.
In particular, if we consider particular $d+1$ points, say, $p_1,\dotsc,p_{d+1}$ as $P$ evolves to $P'$, then the
number of orientation changes of these points is equal to the number of zeros of the polynomial
\[
f(x)\eqdef \det
\begin{bmatrix}
  1&1&\cdots&1\\
  p_{1,1}+x p_{1,1}''&p_{2,1}+x p_{2,1}''&\cdots&p_{d+1,1}+x p_{d+1,1}''\\
  p_{1,2}+x p_{1,2}''&p_{2,2}+x p_{2,2}''&\cdots&p_{d+1,2}+x p_{d+1,2}''\\
  \vdots&\vdots&\ddots&\vdots\\
  p_{1,d}+x p_{1,d}''&p_{2,d}+x p_{2,d}''&\cdots&p_{d+1,d}+x p_{d+1,d}''
\end{bmatrix}\]
on the interval $(0,+\infty)$. Note that the coefficients of the polynomial depend on $\lambda=(\lambda_1,\dotsc,\lambda_d)$.
Substituting $P''=DP'$, we can write the polynomial as
\begin{equation}\label{eq:fdet}
f(x)=\det
\begin{bmatrix}
  1&1&\cdots&1\\
  p_{1,1}+x \lambda_1 p_{1,1}'&p_{2,1}+x \lambda_1 p_{2,1}'&\cdots&p_{d+1,1}+x \lambda_1 p_{d+1,1}'\\
  p_{1,2}+x \lambda_2 p_{1,2}'&p_{2,2}+x \lambda_2 p_{2,2}'&\cdots&p_{d+1,2}+x \lambda_2 p_{d+1,2}'\\
  \vdots&\vdots&\ddots&\vdots\\
  p_{1,d}+x \lambda_d p_{1,d}'&p_{2,d}+x \lambda_d p_{2,d}'&\cdots&p_{d+1,d}+x \lambda_d p_{d+1,d}'
\end{bmatrix}.
\end{equation}

Since the polynomial $f$ is of degree $d$, it has at most $d$ roots. This means that the number of orientation changes
is at most $d$. Since this holds for any $d+1$ points of $P$, it follows that $\dist(P,P')\leq d\binom{n}{d+1}$.

\paragraph{Case $d$ is even:}
This case is particularly easy. Indeed, in this case both $\lambda=(1,1,\dotsc,1)$ and $\lambda=(-1,-1,\dotsc,-1)$
are the valid choices for the scaling factors. The two resulting polynomials, which we call $f^{(1)}$ and $f^{(2)}$,
are related by $f^{(1)}(x)=f^{(2)}(-x)$. Hence, the total number of roots of $f^{(1)}$ and $f^{(2)}$ on the interval $(0,+\infty)$
is equal to the number of roots of either of them on $\R\setminus\{0\}$, which is at most $d$.
Since this holds for all choices of $d+1$ points, it follows that $2\dist(P,P')=\dist(P,P')+\dist(P,-P')\leq d\binom{n}{d+1}$,
proving \Cref{thm:genupper} in this case.

\paragraph{Case $d$ is odd:}
Pick a sequence of real numbers $\lambda_1,\dotsc,\lambda_d$ whose absolute values rapidly decay, i.e. $1\gg \abs{\lambda_1}\gg \dotsb \gg \abs{\lambda_d}$.
We shall choose the signs of these real numbers later.

Because the scaling factors $\lambda_i$ rapidly decay, it is easy to estimate the coefficients of $f$.
\begin{lemma}\label{lem:rj}
Let the coefficients of $f$ be $f(x)=\sum_{j=0}^d c_j x^j$. Then
\[
  c_j=\lambda_1\dotsb\lambda_j r_j(1+O(\lambda_{j+1}/\lambda_j)),
\]
where
\[
r_j\eqdef \det \begin{bmatrix}
  1&1&\cdots&1\\
  p_{1,1}'&p_{2,1}'&\cdots&p_{d+1,1}'\\
  p_{1,2}'&p_{2,2}'&\cdots&p_{d+1,2}'\\
  \vdots&\vdots&\ddots&\vdots\\
  p_{1,j}'&p_{2,j}'&\cdots&p_{d+1,j}'\\
  p_{1,j+1}&p_{2,j+1}&\cdots&p_{d+1,j+1}\\
  \vdots&\vdots&\ddots&\vdots\\
  p_{1,d}&p_{2,d}&\cdots&p_{d+1,d}
\end{bmatrix}.
\]
\end{lemma}
\begin{proof}
  This formula for $r_j$ follows by expanding the determinant in \eqref{eq:fdet} and keeping the largest terms.
\end{proof}
By perturbing $P'$ slightly we can ensure that none of the determinants defining the constants $r_j=r_j(p_1,\dotsc,p_{d+1},p_1',\dotsc,p_{d+1}')$
vanish. Fix such a perturbation.

Let $I_j$ be the real interval between $-\frac{r_{j-1}}{r_j}\cdot \frac{1}{2\lambda_j}$ and $-\frac{r_{j-1}}{r_j}\cdot \frac{2}{\lambda_j}$.
If the $\lambda$'s grow sufficiently quickly, then the polynomial $f$ has a root in $I_j$ for each $j\in [d]$.
Indeed, the value of $f$ on $I_j$ is dominated by the terms $\lambda_1\dotsb\lambda_{j-1}\lambda_j r_jx^j$
and $\lambda_1\dotsb\lambda_{j-1} r_{j-1}x^{j-1}$, and sum of these two terms changes sign in this interval.
Since $f$ has at most $d$ roots and the $d$ intervals are disjoint, these are all the roots of $f$.

Importantly, this implies that the number of orientation changes
in the linear motion of points $p_1,\dotsc,p_{d+1}$ is equal to
the number of negative elements among $\{\lambda_j r_{j-1}r_j  : j\in [d]\}$.

Pick the sign vector of $\lambda_1,\dotsc,\lambda_d$ uniformly at random among $\{-1,+1\}$ vectors
of length $d$ with even number of $-1$'s. Since $d>1$, the sign of any individual $\lambda_i$ is uniform,
and so the expected number of orientation changes is $d/2$. Since this holds for any $d+1$ points,
the bound $\dist(P,P')\leq \frac{d}{2}\binom{n}{d+1}$ follows.

\section*{Proof of \texorpdfstring{\Cref{thm:sametype}}{Theorem 2}}
To prove \Cref{thm:sametype} we use a linear motion as in the proof of \Cref{thm:genupper}, combined with a suitable pre-processing that guarantees certain tuples do not change order type during the linear motion.
We first require a lemma.

\begin{lemma}\label{lem:rotation}
Let $Q = (q_1,\ldots, q_{d+1})\in (\R^d)^{d+1}$ be a tuple consisting of the vertices of a regular simplex inscribed in the unit sphere, and let $\rho\in \SO(d)$ be a rotation. 
If the orientation of the simplex does not remain constant during the linear motion from $Q$ to $\rho(Q)$, then $-\rho$ has a fixed point other than the origin.
\end{lemma}
\begin{proof}
Let $t\in(0,1)$. 
Then $tQ + (1-t)\rho(Q)$ is affinely dependent if and only if there exist $\alpha_1,\ldots, \alpha_{d+1}$ not all zero so that $\sum_{i=1}^{d+1}\alpha_i = 0$, and \[
0 = \sum_{i=1}^{d+1} \alpha_i\Big(tq_i + (1-t)\rho(q_i)\Big).
\]
Let $v = \sum_{i=1}^{d+1} \alpha_i q_i$, and note that the latter equality can be expressed as $tv = (t-1)\rho(v)$.
The condition $\sum_{i=1}^{d+1} \alpha_i=0$ implies that not all $\alpha$'s are equal, and so $v\neq 0$.
Since $\rho$ preserves norms and $tv = (t-1)\rho(v)$, we see that $t=\frac{1}{2}$, and hence $v = -\rho(v)$. 
\end{proof}

We note that the above lemma is interesting only in even dimensions.
Indeed, any rotation in odd-dimensional space fixes a line through the origin by the hairy ball theorem, and so the conclusion of the lemma is trivial when $d$ is odd. 
On the other hand, when $d$ is even we can find $\rho\in \SO(d)$ with neither $\rho$ nor $-\rho$ having a nonzero fixed point. 
One choice is the rotation induced by a cyclic permutation of the vertices of a simplex, but in fact a generic choice suffices: the set of rotations with $\rho$ or $-\rho$ having a nonzero fixed point is a closed subset of $\SO(d)$ with positive codimension. 

\paragraph{Main proof.} Using the above observations, we now prove \Cref{thm:sametype}. 
First, perform a motion so that the first $d+1$ points of $P$ are vertices of a regular simplex inscribed in the unit sphere, noting that this can be accomplished at no cost by rotations and linear scalings that move one point at a time to an appropriate vertex on the sphere. 
Let $Q$ denote these vertices, and choose a rotation $\rho$ so that both $\rho$ and $-\rho$ have no nonzero fixed point. 
As $P'$ has the same order type as $P$, we may arrange that $\rho(Q)$ comprises the first $d+1$ points of $P'$. 

\Cref{lem:rotation} now guarantees that during the linear motions from $P$ to $P'$ and from $P$ to $-P'$ the first $d+1$ points never change orientation, and the same will be true of any sufficiently small perturbation of $P$ and $P'$.
Hence we may proceed similarly to the proof of \Cref{thm:genupper}, counting the zeroes of polynomials associated to all tuples except for the one consisting of the first $d+1$ points.
In this way, we find that $2\dist(P,P') = \dist(P,P') + \dist(P, -P') \le d\left(\binom{n}{d+1} - 1\right)$, which proves the theorem. 

\section*{Proof of \texorpdfstring{\Cref{thm:nonelongated}}{Theorem 3}}
  Suppose $P=(p_1,\dotsc,p_{d+1})$ and $P'=(p_1',\dotsc,p_d')$ are two $(d+1)$-tuples of points
  in general position with the same orientation. Call a rotation $\rho\in \SO(d)$ \emph{good} if for the linear motion from $P$ to $\rho P'$ the polynomial $f$ from \eqref{eq:fdet} has at most  
  $d/2$ zeros on $(0,\infty)$. 
   Since $d\equiv 2\pmod 4$, the number $d/2$ is odd, which implies that, for good $\rho$, the number of zeros of $f$ on $(0,\infty)$ is at most $(d-1)/2$
  since it must be even.

  Let $G(P,P')\eqdef \{\rho \in \SO(d) : \rho\text{ is good}\}$.
  Endow $\SO(d)$ with the natural probability measure, and denote by $m(P,P')$ the measure of $G(P,P')$.
  Recall that a function $f\colon X\to (-\infty,\infty]$ is \emph{lower semicontinuous} if
  the preimage of every interval of the form $(\alpha,\infty]$ is open.
  \begin{lemma}\label{lem:semicont}
    Let $m(P,P')$ be as above. Then $m(P,P')>1/2$, and the function $(P,P')\mapsto m(P,P')$ is lower semicontinuous. 
  \end{lemma}
  \begin{proof}
    Since the roots of a polynomial depend continuously on its coefficients and the interval $(0,\infty)$ is open, the set
    \[
      \mathcal{G}\eqdef \{(P,P',\rho) :\rho\in G(P,P')\}
    \]
    is an open
    subset of $(\R^d)^{d+1}\times (\R^d)^{d+1}\times \SO(d)$. Consequently, the set $\overline{G(P,P')}\eqdef \SO(d)\setminus G(P,P')$
    is closed.

    From the proof of \Cref{thm:genupper} we know that $\rho\in G(P,P')$ or $-\rho\in G(P,P')$, for all $\rho$; 
    hence $\overline{G(P,P')}$ is disjoint from $-\overline{G(P,P')}$.
    Since $\SO(d)$ is connected, $\overline{G(P,P')}\cup(-\overline{G(P,P')})$ is a closed proper subset of $\SO(d)$ and must have measure strictly less than $1$, implying that $m(P,P') >1/2$.

    The lower semicontinuity of $(P,P')\mapsto m(P,P')$ follows from the fact that $\mathcal{G}$ is open
    and the spaces $\SO(d)$ and $(\R^d)^{d+1}\times (\R^d)^{d+1}$ are locally compact Hausdorff spaces, see \cite[Proposition 7.6.5]{cohn_measure}.
  \end{proof}
  Let $\Delta$ be the simplex spanned by the tuple $P=(p_1,\dotsc,p_{d+1})\in (\R^d)^{d+1}$.
  We define the \emph{aspect ratio} of the simplex to be $a(P)\eqdef \diam(\Delta)^d/\vol(\Delta)$.
  The aspect ratio is invariant under rigid motions and scaling.
  \begin{lemma}\label{lem:abound}
    For every $B$, there is a $\delta>0$ such that
    $m(P,P')>1/2+\delta$ whenever both $a(P)$ and $a(P')$ are bounded above by $B$.
  \end{lemma}
  \begin{proof}
    From the formula \eqref{eq:fdet} for $f(x)$ it is clear that $G(P,P')$ (and hence $m(P,P')$)
    is invariant under scaling of $P'$. By symmetry, it follows that it is also invariant under scaling
    of~$P$. So, it suffices to consider simplices contained in, say, the closed unit ball $B(0,1)$.
    The set of simplices inside the ball that have aspect ratio at most $B$ is compact. Since
    a lower semicontinuous function attains its infimum on a compact set, and the function $m$
    strictly exceeds~$1/2$, the result follows.
  \end{proof}

  \begin{lemma}\label{lem:rand}
    If a tuple $P\in (\R^d)^n$ is $\alpha$-non-elongated, and we pick $d+1$ points $q_1,\dotsc,q_{d+1}$ uniformly and independently at random from $P$, then
    $\Pr[a(q_1,\dotsc,q_{d+1})>B]\to 0$ as $B\to \infty$.
  \end{lemma}
  \begin{proof}
    Without loss of generality the minimum distance between points of $P$ is $1$.
    By the usual packing argument, this implies that a ball of radius $\sqrt[d]{n/B}$ contains
    $O(n/B)$ points of $P$. Therefore, $\Pr[\diam(P)\leq \sqrt[d]{n/B}]\leq \Pr[\abs{q_1-q_2}\leq \sqrt[d]{n/B}]\to 0$ as
    $B\to\infty$.

    Sample the points $q_1,\dotsc,q_d$ first, and look at the hyperplane $H$ that they span. Since all pairwise distance are at most
    $\alpha\sqrt[d]{n}$, the volume of the simplex $\conv(q_1,\dots, q_{d+1})$ is less than $n/B$ only if the point $q_{d+1}$ lies in a slab 
    of width $O(\alpha^{d-1}n/ B)$ around $H$. The probability of this event is at most
    $O(\alpha^d n/8^d B)$, and in particular tends to zero as well. So, with high probability both $\diam(P)\leq \sqrt[d]{n/B}$
    and $\vol(\conv(q_1,\dots, q_{d+1}))\geq n/B$; in particular, $a(P)\leq B^2$.
  \end{proof}

  \Cref{thm:nonelongated} now easily follows. Select $B$ sufficiently large so that the probability in \Cref{lem:rand}
  is at most $1/4$. Let $X$ be the set of all $(d+1)$-tuples whose aspect ratio both in $P$ and $P'$
  is at most $B$. We have $\abs{X}\geq \tfrac{1}{2}\binom{n}{d+1}$.
  Let $\delta$ be the constant from \Cref{lem:abound} for this value of $B$.
  Pick $\rho$ be a uniform random element of $\SO(d)$, and consider the linear motion from $P$ to $\rho P'$.
  By \Cref{lem:abound}, the expected cost of this motion on a tuple in $X$ is at most $\frac{d}{2} - \delta$. 
  Therefore the expected total cost of this motion is at most
  \[
  \frac{d}{2}\left(\binom{n}{d+1}-\abs{X}\right) + \left(\frac{d}{2}-\delta\right)\abs{X} \,\,=\,\,    \frac{d}{2}\binom{n}{d+1}-\delta\abs{X}\,\,= \,\,\left(\frac{d}{2} - \frac{1}{2}\delta\right) \binom{n}{d+1} .
  \]

  \paragraph{Acknowledgment.} The proof of \Cref{lem:semicont} would not be as short without the aid of Gautam Iyer who
  directed us to \cite{cohn_measure}.

\bibliographystyle{plain}
\bibliography{ordertypes.bib}
\end{document}